\documentclass[12pt,reqno]{amsart}

\usepackage{amsmath, amssymb}
\usepackage[latin1]{inputenc}

\newtheorem{theorem}{Theorem}[section]
\newtheorem{lemma}[theorem]{Lemma}
\newtheorem{prop}[theorem]{Proposition}
\newtheorem{cor}[theorem]{Corollary}
\newtheorem{rmk}[theorem]{Remark}

\newcounter{defn}

\newcommand{\sect}{\vspace{3mm} \setcounter{equation}{0} \setcounter{defn}{0} \section}

\newcommand{\w}[1]{\langle {#1} \rangle}
\newcommand{\pf}{\noindent {\bf Proof. \hspace{2mm}}}
\newcommand{\ef}{ \hfill $ \blacksquare $ \vskip 3mm}
\newcommand{\be}{\begin{equation}}
\newcommand{\ee}{\end{equation}}
\newcommand{\bea}{\begin{eqnarray}}
\newcommand{\eea}{\end{eqnarray}}

\newcommand{\Bsigma}{\bf \Sigma}

\newcommand{\bC}{{\mathbb C}}
\newcommand{\bR}{{\mathbb R}}
\newcommand{\bN}{{\mathbb N}}

\newcommand{\bS}{{\bf S}}

\newcommand{\vF}{{\mathcal F}}

\newcommand{\vL}{{\mathcal L}}

\newcommand{\vN}{{\mathcal N}}

\newcommand{\vW}{{\mathcal W}}

\newcommand{\vS}{{\mathcal S}}

\begin{document}

\title[Number of eigenvalues]{Number of eigenvalues for a class of non-selfadjoint Schr\"odinger operators}
\author{ Xue Ping WANG}
\address{Département de Mathématiques\\
UMR 6629 CNRS\\
Université de Nantes \\
44322 Nantes Cedex 3  France \\
E-mail: xue-ping.wang@univ-nantes.fr
}
\date{}
\thanks{Research  supported in part by the French National Research Project NONAa, No. ANR-08-BLAN-0228-01, on {\it Spectral and microlocal analysis of non-selfadjoint operators}.}

\begin{abstract}
In this article, we prove the finiteness of the number of eigenvalues for a class of Schr\"odinger operators $H = -\Delta + V(x)$  with a complex-valued potential $V(x)$ on $\bR^n$, $n \ge 2$.  If  $\Im V$ is sufficiently small, $\Im V \le 0$ and $\Im V \neq 0$, we show that  $N(V) = N( \Re V )+ k$,  where $k$ is the multiplicity of the zero resonance of the selfadjoint operator $-\Delta + \Re V$ and $N(W)$ the number of eigenvalues of $-\Delta + W$, counted according to their algebraic multiplicity.
\end{abstract}

\maketitle

\sect{Introduction}

Consider the Schr\"odinger operator $H = -\Delta + V(x)$ with a complex-valued potential $V(x) = V_1(x) - iV_2(x)$  on $L^2(\bR^n)$,   where $V_1$ and $ V_2$ are real measurable functions.  $V$  and $H$ will be called dissipative if $V_2(x) \ge 0$ and $V_2(x) >0 $ on some non trivial open set. Assume  that $V$ is a $-\Delta$-compact perturbation.   $H$  is then closed with domain $D(H) = D(-\Delta)$. Let $\sigma(H)$ (resp., $\sigma_{ess}(H)$, $\sigma_{pp}(H)$) denote the spectrum (resp., essential spectrum, eigenvalues) of $H$. By  Weyl's  essential spectrum theorem and the analytic Fredholm theorem, one has $\sigma_{ess}(H)  = [0, \infty[$ and the spectrum of $H$ is discrete in $\bC\setminus [0, \infty[$, consisting of eigenvalues with finite algebraic multiplicity.  {\it A priori}, the complex eigenvalues of $H$ may accumulate to any non negative real number. \\

For real-valued potentials $V$, it is well-known that the number of the eigenvalues of $H=-\Delta + V(x)$ is finite if $V(x)$ decays like $O(|x|^{-\rho})$ for some $\rho>2$  and some universal estimates hold. See \cite{rs} for a recent survey on this topic. The study of the number of eigenvalues for Schr\"odinger operators with complex-valued potentials has also a long history (\cite{dav,nai,sch}). But there are relatively fewer results and most of them are concerned with one-dimensional problems.  In particular, B. S. Pavlov studied  in \cite{pav1,pav2}   the Schr\"odinger operator with a complex-valued potential on the positive half-line under the non-selfadjoint boundary condition $y'(0) -h y(0) =0$ for some complex parameter $h$. He proved that the number of the eigenvalues is finite  if  $|V(x)|= O(\exp(-\varepsilon|x|^{1/2}))$ for some $\varepsilon>0$, as $x\to +\infty$,  and that for each $\lambda>0$, there exist some value of $h$ and some real smooth function $V(x)$ decaying like $O (\exp(-c|x|^{\beta}))$  at the infinity for some $c>0$ and $\beta <1/2$ such that  the eigenvalues  accumulate to the point $\lambda$.  Pavlov's ideas  are utilized in \cite{tb} to prove the finiteness of the number of eigenvalues for non-selfadjoint Schr\"odinger operators on the whole real axis  and in \cite{ge}  for the limit set of eigenvalues of complex Jacobi matrices.  In \cite{aad}, the authors showed that for potentials $V$ satisfying $\|V(x) e^{\gamma |x|} \|_1 <\infty$ for some $\gamma>1$ sufficiently large, all the eigenvalues satisfy the estimate $|\lambda| \le \frac 9 4 \|V\|_1^2$. In multidimensional case, there are several recent works on the estimate of eigenvalues for non-selfadjoint Schr\"odinger operators. In \cite{bo,flls}, Lieb-Thierring type estimates for Schr\"odinger operators  are  proved for eigenvalues lying outside the sector $\{z; |\Im z| < t \Re z\}$, $t>0$,  with a $t$-dependent constant. In  \cite{ls}, A. Laptev and O. Safronov obtained several estimates on some sums of the imaginary parts of eigenvalues. In particular, their results give some rate of convergence of the imaginary parts of eigenvalues and imply that if $n=3$ and $V_2 \ge 0 $ is integrable, the eigenvalues of $-\Delta - i V_2$ can not accumulate to zero.
In this work, we will prove that for a class of complex-valued potentials on $\bR^n$ with $n \ge 2$, the eigenvalues of $H = -\Delta + V(x)$ can not accumulate to a non negative real number and that if $V$ is dissipative,  the zero eigenvalue and the  zero resonance (in the sense of \cite{jk}) of the selfadjoint operator $ H_1= -\Delta + V_1$  are turned into complex eigenvalues of $H = -\Delta + V_1 -i V_2$, once a small dissipative part of the potential, $-iV_2$, is added on.

The minimal assumptions used in this work are the following. $V_1$ and $V_2$  satisfy the estimates
\begin{equation} \label{ass1}
\mbox{ $V_j$ is $-\Delta$-compact and $ |V_j(x)| \le C \w{x}^{-\rho_j}$ for $|x| > R $,}
\end{equation}
for some $R>0$ and $\rho_j>1$, $j=1,2$. Here $\w{x} = (1 + |x|^2)^{1/2}$.
For the reason of convenience, we assume also that
$V_j \in L_{\rm loc}^p(\bR^n)$ with $p = \frac n 2$ if $ n \ge 3$, $p>1$ if $n=2$ and $p=2$ if $n=1$. This allows us to apply the known results on the unique continuation theorem  to $H=-\Delta + V$.
The results on the finiteness of eigenvalues will  be proved  for $ n\ge 2$ and under stronger conditions on $V$ with $\rho_2>2$ and the real part of the potential,  $V_1$, satisfying
\begin{equation} \label{ass3}
V_1 (x) = \frac{q(\theta)}{r^2} + O(\w{x}^{-\rho_1'}), \quad |x| > R,
\end{equation}
for some $R>0$ and  $\rho_1'>2$, where $r =|x|$,  $x = r\theta$  is spherical coordinates on $\bR^n$, $n \ge 2$, and $q(\theta)$ is a real continuous function on $\bS^{n-1}$ such that the lowest eigenvalue, $\mu_1$,  of $-\Delta_{\bS^{n-1}} + q(\theta)$ on $\bS^{n-1}$ verifies
\begin{equation} \label{ass4}
\mu_1 > - \frac{(n-2)^2}{4}.
\end{equation}
Note that if $n\ge 3$ and $V_1$ satisfies (\ref{ass1}) for some $\rho_1>2$, (\ref{ass3}) and (\ref{ass4}) are satisfied with $q=0$ and $\rho_1' =\rho_1$. For $n=2$,  the condition (\ref{ass4})  requires the potential to be positive in some sense when $r$ is large enough and potentials with compact support are excluded.  For potentials satisfying (\ref{ass1}), (\ref{ass3}) and (\ref{ass4}), we shall say that zero is a regular point of $H$, if for any $u \in H^{1,-s}$ (the  weighted  first order Sobolev space with the weight $\w{x}^{-s}$) for some $s<1$  such that $Hu=0$, one has $u=0$. We shall show in Lemma \ref{lem2.2} that if $V$ is dissipative, zero is a regular point of $H$.  A solution $u$ of $Hu=0$ is called a resonant state if  $u \in H^{1,-s}\setminus L^2$ for some appropriate $s>0$. See \cite{jk}.

\begin{theorem} \label{th1.1} (a). Let $n \ge 1$. Under the condition (\ref{ass1}), the eigenvalues of $H$ can not accumulate to a positive real number $E>0$.

(b). Let $ n \ge 2$. Under the conditions (\ref{ass1})-(\ref{ass4}) with $\rho_1'>2$ and $\rho_2>2$,  suppose in addition that zero is a regular point of $H$. Then, the number of the eigenvalues of $H$ is finite. In particular, if $V$ is dissipative and (\ref{ass1})-(\ref{ass4}) are satisfied with $\rho'_1>2$ and $\rho_2>2$, the number of the eigenvalues of $H$ is finite.
\end{theorem}

Part (a) of Theorem \ref{th1.1} says that for Schr\"odinger operators with a short-range complex-valued potentials, zero is the only possible point of accumulation of the eigenvalues. Part (b) of Theorem \ref{th1.1}  is optimal in the sense that even in the case $V_2=0$, if $V_1$ is of critical decay at the infinity, but the condition (\ref{ass4}) is not satisfied,   $H_1= -\Delta + V_1$ may have an infinite number of negative eigenvalues accumulating to zero. 

The next objective of this work is to estimate the number of eigenvalues of a non-selfadjoint Schr\"odinger operator  when the imaginary part of the potential is small.  Denote $H(\lambda) = H_1 - i \lambda V_2$ where $\lambda \in \bR $ is a small parameter. Let $N(\lambda) $ (resp. $N_1$) be the total number of the complex eigenvalues of $H(\lambda)$ (resp., $H_1$), counted according to their algebraic multiplicity. It is easy to show that under the same conditions, if $0$ is a regular point of $H_1$, then
 \begin{equation} \label{NV}
 N(\lambda) = N_1
 \end{equation}
for $ |\lambda |\le \lambda_0$.  See Section 4. A more interesting question is the case when zero happens to be an eigenvalue and a resonance of $H_1$.  For critical potentials of the form (\ref{ass3}), the zero resonance may appear in any space dimension $n$ with arbitrary multiplicity depending on the small eigenvalues of $-\Delta_{\bS^{n-1}}+ q$ (\cite{w2}).  In this work, we will only study the particular case  when $n \ge 3$ and
\begin{equation}\label{mu1}
\left\{\nu =\sqrt{ \mu + \frac{(n-2)^2}{4}}, \mu \in \sigma ( - \Delta_{\bS^{n-1}} + q) \right\} \cap ]0, 1] =\left\{ \nu_1  \right \},
\end{equation}
where $\nu_1 =  \sqrt{ \mu_1 + \frac{(n-2)^2}{4}}$.
The condition (\ref{mu1}) ensures, among others, that if zero is a resonance of $H_1$, then it is simple.
If
 \[
 \left\{\nu =\sqrt{ \mu + \frac{(n-2)^2}{4}}, \mu \in \sigma ( - \Delta_{\bS^{n-1}} + q) \right\} \cap ]0, 1] =\emptyset,
 \]
 zero is never a resonance of $H_1$. Let $\varphi_0$ be a normalized eigenfunction of $- \Delta_{\bS^{n-1}} + q$ associated with $\mu_1$, which can be taken to be positive. Set
\[
W_1(x) = V_1 (x) - \frac{q(\theta)}{r^2}, \quad \eta_0(x) = \frac{\varphi_0(\theta)}{r^{\frac{n-2}{2} - \nu_1}}, \quad x =r\theta.
\]

\begin{theorem} \label{th1.2}  Let $n \ge 3$. Assume (\ref{ass1}) - (\ref{ass4}) with $\rho_1' >3$ and $\rho_2>3$ and (\ref{mu1}).  Suppose in addition that $V$ is dissipative: $V_2 \ge 0$ and $V_2 \neq 0$. Let zero be an eigenvalue \mbox{\rm or} a resonance of $H_1$. When zero is a resonance, suppose that
 \begin{equation} \label{ass5}
\nu_1 \in [1/2, 1]  \quad \mbox{ and }  \quad \overline{ \w{W_1 \eta_0 , \phi} } \w{V_2 \eta_0, \phi} <0,
 \end{equation}
where $\phi$ is a resonant state. Then one has
 \begin{equation}
 N(\lambda) =  N_1 +k,
 \end{equation}
 for $0< \lambda <\lambda_0$. Here $k$ is equal to $0$ if zero is not a resonance of $H_1$, and is equal to $1$ if zero is a resonance of $H_1$ and $N_1$ is the total number of eigenvalues of $H_1$, including the zero ones.
\end{theorem}

Note that $ \w{W_1 \eta_0 , \phi} \neq 0$ if $\phi$ is a resonant state (\cite{w2}) and the condition (\ref{ass5}) is independent of the choice of $\phi$. (\ref{mu1}) is always satisfied if $q $ is an appropriate constant. In particular, if $n=3$ or $4$ and the condition  (\ref{ass1})  is satisfied with $\rho_j >2$, $j=1,2$,  one has   $q =0$, $\nu_1 = \frac{1}{2}$ or $1$, respectively,
$W_1 = V_1$ and $\eta_0 $ is constant: $\eta_0 = \frac{1}{\sqrt{|\bS^{n-1}|}}$. The condition (\ref{ass5}) is simplified  as
\begin{equation}
\overline{ \w{V_1, \phi} } \w{V_2, \phi} <0,
\end{equation}
Therefore, a particular case of Theorem \ref{th1.2} is the following

\begin{cor} \label{cor1.3} Let $n = 3, 4$ and  let $V_1$ verify the condition (\ref{ass1}) with $\rho_1>3$ and $V_1 \le 0, V_1 \neq 0$. Let $V_2 = -V_1$ and $H(\lambda) = -\Delta + (1+ i \lambda) V_1$, where $\lambda $ is a small positive parameter. Assume that zero is an eigenvalue or a resonance of $H_1$. Then one has
\begin{equation}\label{exemple}
N(\lambda) =\left\{ \begin{array}{cc}
N_1, & \mbox{ if zero is not a resonance of $H_1$};\\
N_1 + 1, & \mbox{ if zero is  a resonance of $H_1$}.
\end{array}
\right.
\end{equation}
for $\lambda>0$ sufficiently small.
\end{cor}

In fact, if zero is  a resonance, one has
\[
\overline{ \w{V_1, \phi} } \w{V_2, \phi}  = -| \w{V_1, \phi} |^2 <0
\]
for any resonant state $\phi$, because $\w{V_1, \phi} \neq 0$. (\ref{ass5}) is satisfied. An example of the potential $V_1$ for which zero is a resonance but not an eigenvalue of $H_1$ is given in Section 4.

The plan of this work is as follows. In Section 2, we prove that for short-range complex-valued potential  $V$,  the eigenvalues of $H = - \Delta + V(x)$ can not accumulate to a positive real number $E>0$.  In  Section 3,  we show that for the class of potentials satisfying the conditions of Theorem \ref{th1.1} (b), the eigenvalues of $H$ can not accumulate to zero. The number of eigenvalues under perturbations is discussed in Section 4. The main attention is payed to the case where $V$ is dissipative and  the zero eigenvalue and the zero resonance of $H_1$ are present.
Theorem \ref{th1.2}  follows as a consequence of this analysis.
\\

{\noindent \bf Notation.} $H^{r,s}$, $ r, s
\in \mathbb{R}  $, denotes the weighted Sobolev space of order $r$
defined by $H^{r, s} =\{ f \in \vS'(\bR^n); \w{x}^{s} (1- \Delta)^{r/2}f \in L^2\}$ equipped with the natural  norm noted as $\|\cdot\|_{r,s}$. The dual product between
$H^{r,s}$ and $H^{-r,-s}$ is identified with $L^2$-scalar product. Denote
$H^{0,s} = L^{2,s}$ and $H^{r,0} = H^r$.  $\vL(r,s; r', s')$ is the space of continuous linear operators from $H^{r,s}$ to $H^{r^{\prime},s^{\prime}}$ and $\vL(r,s) = \vL(r,s; r, s)$.
\\

\noindent {\bf Acknowledgements.}  {\it  The author should like to thank V. Bruneau and J. Sj\"ostrand for useful discussions on this subject, and J.-F. Bony and D. Yafaev for bringing to his attention the references \cite{aad} and \cite{pav1}, \cite{pav2}, respectively.}

\sect{Absence of eigenvalues near a positive real number}

The condition (\ref{ass1}) implies that $V_j$, $j=1,2$, are bounded as operators from $H^{1, s}$ to $H^{-1, s +\rho_j}$. It will be convenient to regard $H$ as a bounded operator from $H^{1, s}$ to $H^{-1, s}$ for any $s\in \bR$.  The condition (\ref{ass1}) with $\rho_j>1$ for $V_j$ guarantees, among others, that the positive eigenvalues of $H$ are absent.  Let $R_0(z) = (-\Delta -z)^{-1}$. It is well-known  that the boundary values of the resolvent
\begin{equation}
R_0(\lambda \pm i0) =\lim_{\epsilon \to 0_+} R_0(\lambda \pm i \epsilon), 
\end{equation}
exist in $\vL(-1, s; 1,-s)$ for any $s>1/2$ and is continuous in $\lambda >0$.

\begin{lemma} \label{lem2.1}
Assume (\ref{ass1}).
Let $z_0 \in \bC$ and $u \in H^{1,-s}$ for some $s>0$  such that $ Hu = z_0u$.

(a). For $z_0= \lambda  > 0$, assume that  $1/2 < s < \min\{ \rho_1/2, \rho_2/2\}$ and that $u$ verifies one of the radiation conditions
\begin{equation}\label{rad}
u = - R_0( \lambda \pm i0) (Vu).
\end{equation}
Then $u=0$.

(b). If $\Im z _0 \neq 0$, then $u$ is rapidly decreasing: $u\in H^{2, r}$ for any $r\ge 0$.
\end{lemma}
\pf  (a). Note that the condition on $V$ and the equation $Hu = z_0 u$  with $u\in H^{1,-s}$ imply that $u\in H^{2, -s}$.  Let $z_0 = \lambda>0$. To fix the idea, suppose that $u$ verifies the incoming radiation condition:
\[
u =  - R_0( \lambda - i0) (Vu).
\]
Since $\rho_j >2s>1$, $u$ can be decomposed as 
\begin{eqnarray}
u& =& v  + w_1, \quad  \nabla u  =   -i \theta \sqrt{\lambda} v   + w_2 \mbox{ with }\\
v & = & \frac{  a_0(\theta)}{ r^{\frac{n-1}{2}}}  e^{-i \sqrt{\lambda} r }
\end{eqnarray}
 for some $ a_0 \in L^2(\bS^{n-1})$ and $w_j \in L^2_{\rm loc}(\bR^n)$  satisfying for some $\epsilon_0>0$ and $R_0>0$
 \[
 |w_j(x)| \le C \w{x}^{- (n-1)/2 - \epsilon_0},
 \]
  for  $|x| >R_0$.  See \cite{is}. Here $ x = r\theta$ with $r = |x| $, $\theta \in \bS^{n-1}$,  with the convention that $\theta = \mbox{ \rm sgn } x$ if $n=1$.  For $R>1$, set
\[
\chi_R(x) = e^{ -\w{x}/R} \quad \mbox{ and } \quad u_R = \chi_R u.
\]
Then, $u_R \in H^{2, \tau}$ for any $\tau >0$. By the formula (14.7.1) of \cite{hor} and an argument of density, one has
\begin{equation}  \label{A}
 4 \lambda \tau \| |x|^\tau u_R\|^2  \le \| |x|^{1 +\tau} (-\Delta -\lambda) u_R\|^2
\end{equation}
for any $\tau>1$ and $R>1$. One can evaluate:
\begin{eqnarray*}
\| |x|^\tau \chi_R v\|^2 & = & \|a_0\|^2  \int_0^\infty r^{2\tau} e^{-2 \w{r}/R} \; dr \\
 & \ge & c_0^2 R^{2 \tau +1} \|a_0\|^2
\end{eqnarray*}
for any $R>1$, where
\[
  c_0 = \left\{ \int_0^\infty r^{2\tau} e^{-2 \w{r}} dr \right\}^{1/2}.
\]
Similarly,  since $w_1$ is in $L^2_{\rm loc}$ and $|w_1(x)|\le C\w{x}^{-(n-1)/2 -\epsilon_0}$ for $|x| >R_0$, one has
\begin{eqnarray*}
\| |x|^\tau \chi_R w_1 \| & \le & R_0^{\tau} \| w_1\|_{L^2( |x| <R_0)} + C \left\{ \int_{R_0}^\infty r^{2(\tau -\epsilon_0)} e^{-2 \w{r}/R} \; dr \right \}^{1/2}\\
 & \le & C_1 (1+  R^{(\tau-\epsilon_0) + 1/2} )
\end{eqnarray*}
for some $C_1>0$ independent of $\tau$ and $R$. This shows that there exists  $R_1>1$ such that
\begin{equation} \label{B}
\| |x|^\tau u_R \| \ge    R^{ \tau +1/2}  c_0/2 ( \|a_0\| - R^{-\epsilon_0/2})
\end{equation}
for all $\tau>1$ and $R>R_1$.
For the right hand side of (\ref{A}), one uses the expression
\[
(-\Delta -\lambda) u_R = - V u_R - (2 \nabla \chi_R \cdot \nabla + \Delta \chi_R) u
\]
Since $\Delta \chi_R = O(R^{-2}) \chi_R$ and $V (x) = O(\w{x}^{-\rho})$ for some $\rho>1$,  a similar calculation gives
\[
\| |x|^{1+\tau} (\Delta \chi_R ) u\| \le C R^{\tau}, \quad \| |x|^{1+\tau} V u_R\| \le C R^{\tau -\rho +3/2}.
\]
Note that
\[
\nabla\chi_R\cdot\nabla = \chi'_R (r) \frac{\partial}{\partial r}.
\]
The radiation condition of $u$ gives
\[
- \nabla \chi_R \cdot \nabla u = \frac{ r}{ R \w{r}} \chi_R (- i \sqrt{\lambda} v + w_2).
\]
Since $\| |x|^{1+\tau} \chi_R v\| \le C R^{(1+\tau) + 1/2}$,  we obtain from the decay of $w_2$ that for some $C, C'>0$,
\begin{equation}
\| |x|^{1+\tau}  \nabla \chi_R \cdot \nabla u \| \le C \sqrt{\lambda} R^{\tau + 1/2} + C' R^{\tau -\epsilon_0 + 1/2}.
\end{equation}
This proves that there exists $C_2, C_3>0$ such that
\begin{equation} \label{C}
\| |x|^{1+\tau}  (-\Delta -\lambda) u_R  \| \le C_2 \sqrt{\lambda} R^{\tau + 1/2} + C_3 R^{\tau -\epsilon'_0 + 1/2}.
\end{equation}
uniformly in $R>R_1$ and $\tau>1$, where $\epsilon_0' = \min\{ \epsilon_0, \rho-1\}$. (\ref{A}), (\ref{B}) and (\ref{C}) show that
\begin{equation}
\sqrt{\lambda \tau} c_0 (\|a_0\| - R^{-\epsilon_0/2}) \le   C_2 \sqrt{\lambda} +  C_3   R^{ -\epsilon'_0}
\end{equation}
for all $\tau>1$ and $R>R_1$. Since $\lambda>0$ and $c_0>0$, this is impossible for $\tau$ and $ R$ sufficiently large, unless $a_0 =0$. Therefore,  $u \in L^{2, -1/2 +\epsilon_0}$ for some $\epsilon_0>0$. This implies that the restriction of the Fourier transform $Vu$ on the sphere $\{ |\xi|^2 =\lambda\}$ is zero and that $u$ is both incoming and outgoing:
\begin{equation}
u =  - R_0(\lambda-i0) Vu = - R_0(\lambda + i0) Vu
\end{equation}
Since $V$ is short-range, the bootstrap argument of \cite{ag} or the microlocal resolvent estimates for $R_0(\lambda \pm i0)$ allow to show that $u \in H^{2, r}$ for any $r>0$. The unique continuation theorem  (\cite{hor,ij}) for the Schr\"odinger equation  implies that   $u=0$. \\

 When $ \Im z_0 \neq 0$,  $(H_1 -z_0)^{-1}$ is continuous from $H^{-1, r}$ to $H^{1, r}$ for any $r\in\bR$. (\ref{ass1}) shows that  $V_2u \in H^{-1, -s+ \rho_2}$ which implies  $ u = i(H_1-z)^{-1}V_2u \in H^{1,-s+ \rho_2}$. Since $\rho_2>0$, an iteration of this argument shows that $u \in H^{1,r}$ for any $r\ge 0$.
\ef

The set of positive real numbers $\lambda >0$ for which  the equation $H u = \lambda u$ admits a non trivial solution in $L^{2,-s}$ for any $s>1/2$ and verifying some radiation condition has been analyzed in the 1970's in relation with the limiting absorption principle for non-selfadjoint Schr\"odinger operators. In \cite{sa,sa2}, it was proved for complex-valued potentials with a long-range real part that this set is bounded with measure zero.  Lemma \ref{lem2.1} shows that for short-range potentials, this singular set is absent.  The following result shows that when $V$ is dissipative,  zero is a regular point of $H$, \emph{i.e., } it is neither an eigenvalue nor a resonance of $H$.

 \begin{lemma} \label{lem2.2} Let $n \ge 2$. Under the condition (\ref{ass1}), suppose in addition that $V_2\ge 0$ and $V_2 \neq 0$.  Let $u \in H^{1,-s}$,  $s < \min\{1, \rho_1/2, \rho_2/2\}$, such that $Hu=0$. Then $u=0$.
\end{lemma}
 \pf
Since $\rho_j \ge 2s$, one has  $ -\Delta u = -Vu \in L^{2, s}$, $ H_1  u = i V_2 u\in L^{2, s}$ and
\begin{equation}\label{real}
\w{ u,  H_1  u} = i \w{ u,  V_2u }.
\end{equation}
We want to show that $\w{H_1 u, u}$ is a real number, although $u$ is not in the domain of the selfadjoint operator $H_1$.
The Fourier transform $\vF$ is a bounded map from $L^{2,r}$ to $H^{r}$ for any $r\in \bR$ and for any $f \in L^{2,r}$, $g \in L^{2, -r}$ one has the Plancherel's formula:
$\w{f, g} = \frac{1}{(2\pi)^n}  \w{\hat{f}, \hat{g}}$, where $\hat{f} $ is the Fourier transform of $f$.  Since the Fourier transform of $-\Delta  u$ is  $ |\xi|^2 \hat{u} $, one has $\hat{u} \in H^{-s}$ and $ |\xi|^2 \hat{u} \in H^{s}$.  It follows that
\[
\w{ -\Delta u, u} = \frac{1}{(2\pi)^n} \w{|\xi|^2 \hat{u}, \hat{u}}.
\]
We are led to verify that $\w{|\xi|^2 \hat{u}, \hat{u}}$ is a real number. It suffices to show that $|\xi| \hat{u} \in L^2$.  Using the characterization of fractional Sobolev norm on $H^s$ (for $0<s<1$) and the generalized Hardy inequality
\[
\| \frac{1}{|x|^\sigma} f\| \le C_\sigma \| f\|_{\dot{H}^\sigma(\bR^n)}, \quad \sigma \in ]0, n/2[,
\]
where $C_\sigma>0$ and $\dot{H}^\sigma$ is the homogeneous Sobolev space,  one can show that
\begin{equation} \label{Du}
\|\w{D}^s |\xi|^{i\mu} \w{D}^{-s}\| \le C_s (1 + |\mu|), \quad   \mu \in \bR,
\end{equation}
where $\w{D} = ( 1 - \Delta_\xi)^{1/2}$. Remark that if $n \ge 3$, this estimate follows easily from 
the standard Hardy inequality and still holds for $s=1$.
Consider the function 
\[
F(z) =\frac{1}{1 +z} \w{f, \w{D}^{s(2z-1)}|\xi|^{2z} \hat{u}}, \quad  z\in \bC, \Re z \in [0, 1],
\]
where $f \in C_0^\infty(\bR^n)$. It follows from (\ref{Du}) that there exists $C>0$ such that
\[
| F( i \mu) |  \le C \| f\| \| \w{D}^{-s} u\|, \quad  | F( 1+ i \mu) |  \le C \| f\| \| \w{D}^{s} u\|,
\] 
for any $f \in  C_0^\infty(\bR^n)$ and $\mu \in \bR$. By the method of complex interpolation, one obtains that
\[
| \w{f, |\xi| \hat{u}} | \le \frac{3C}{2}   \| f\| \| \w{D}^{-s} u\|^{1/2}  \| \w{D}^{s} u\|^{1/2}
\]
for any $f \in  C_0^\infty(\bR^n)$. This proves that $|\xi| \hat{u}$ is in $ L^2$. Therefore,
\[
\w{|\xi|^2 \hat{u}, \hat{u}} = \| |\xi| \hat{u} \|^2  \ge 0.
\]
This proves that $\w{u, H_1u} =  \w{  u,   -\Delta u} + \w{ u, V_1u}$ is a real number.
 It follows from (\ref{real}) that $\w{H_1u, u} =0$ and $\w{V_2 u, u} =0$. Since $V_2 \ge 0$ and $V_2 \neq 0$,  one has $V_2u=0$ and $u(x)=0$ for $x$ in a non trivial open set $\Omega$. Now that $u$ verifies the equation $H_1 u =0$ and $u=0$ on
$\Omega$, we can apply the unique continuation theorem  (\cite{ij}) to $H_1$ to conclude $u=0$ on $\bR^n$. \ef

\begin{rmk}
Let $s<1$ and $W: H^{1,-s} \to H^{-1,s}$ be a bounded operator. The proof of Lemma \ref{lem2.2} shows that if  $\vN$ is a subspace of $H^{1,-s}$ such that   $ (-\Delta + W)\phi=0$ for any $\phi \in \vN$, the form
\begin{equation}
Q(\phi) = \w{ \phi, - \Delta \phi} = \w{ \phi, -W\phi}
\end{equation}
is positive definite on $\vN$, because constant functions do not belong to $H^{1, -s}(\bR^n)$ with $s<1$ when $n\ge 2$.
\end{rmk}

Denote $D(z_0, r) = \{z \in \bC;  |z-z_0| < r  \}$ for $z_0 \in \bC$ and $r>0$,  $D_\pm (z_0, r) = D(z_0, r) \cap \bC_\pm$ and $D'(0,r) = D(0,r)\setminus [0, r[$. \\

\noindent{\bf Proof of Theorem \ref{th1.1} (a).}  Let $E_0 >0$. Assume (\ref{ass1}) with $\rho_j >1$, $j=1,2$. We want to prove that  there exists $\delta>0$ such that
 \begin{equation}
\sigma_{\rm pp}(H) \cap  D(E_0, \delta) =\emptyset.
\end{equation}
Since $H$ has no positive eigenvalues, it suffices to show that $\sigma_{\rm pp}(H) \cap  D_\pm (E_0, \delta) =\emptyset$. We only show that $\sigma_{\rm pp}(H) \cap  D_-(E_0, \delta) =\emptyset$. The other case can be proved in the same way.  Fix $1/2 < s< \rho_j/2$, $j=1,2$. Set $F(z) =  R_0(z)V$ for $z \in \bC_-$ with $\Re z$ near $E_0$. Then $F(z)$ is compact as operators from $H^{1,-s}$ to $H^{1,-s}$  and $F(\lambda)=  R_0(\lambda -i0)V$ exists for $\lambda >0$ real and is Hölder-continuous. By lemma \ref{lem2.1}, $-1$ is not an eigenvalue of $F(E_0)$ in $H^{1, -s}$. The Fredholm theorem shows that $1+ F(E_0)$ is invertible on $H^{1,-s}$.
 Since $F(z) $ is continuous in $ \overline{D}_-(E_0, r_0)$,
there exists some $\delta>0$ such that
\begin{equation}
\| F(z) - F(E_0)\|  \| (1+ F(E_0))^{-1}\| \le c <1, \quad z \in \overline{D}_-(E_0, \delta).
\end{equation}
It follows  that $( 1+ F(z))^{-1} $ exists and
\[
\| (1+ F(z))^{-1} \| \le \frac{\| (1+ F(E_0))^{-1}\|}{1- c}, \quad z \in \overline{D}_-(E_0, \delta).
\]
In particular, $ 1+ F(z)$ is injective on $H^{1,-s}$.  Since $ R_0(z)(H-z) = 1+ F(z)$,  $H-z$ is injective on $H^{1,-s}$. Lemma \ref{lem2.1} (b) implies that $\sigma_{pp}(H)\cap {D}_-(E_0, \delta) =\emptyset$  and  $D_-(E_0, \delta)$ is contained in the resolvent set of $H$.
\ef

From the equation $R(z) = (1+ F(z))^{-1}R_0(z)$ where $R(z) =(H-z)^{-1}$,  it follows from the proof  of Theorem \ref{th1.1} (a) that the boundary values $R(\lambda \pm i0)$ of the resolvent exist and are H\"older-continuous in  $\vL(-1, s; 1,-s)$,  $s>1/2$, for $\lambda >0$. 
In \cite{roy},  the Mourre's theory is generalized to a  class of abstract dissipative operators and semiclassical resolvent estimates are deduced in the presence of trapped trajectories.
In \cite{sa,sa2}, it is proved under more general conditions that the limiting absorption principle holds outside some set of zero measure.

\sect{Finiteness of the number of eigenvalues}

\noindent{ \bf Proof of Theorem \ref{th1.1} (b).}
Using the high energy resolvent estimates for the Laplacian
\[
\|\w{x}^{-s} R_0(z) \w{x}^{-s} \| \le C_s |z|^{-1/2},
\]
for $|z| >1$ and $z\not\in \bR$, for any $s>1/2$, it is easy to show that under the condition (\ref{ass1})  with $\rho_j>1$, there exists $R>0$ such that
\[
\|\w{x}^{-s} R_0(z) V \w{x}^{s}  \| <1
\]
for $1/2 <s < \rho_j/2$ and $z$ with $|z|>R$ and $z\not \in \bR_+$. This means that $1 + R_0(z)V$ is invertible on $L^{2,-s}$.  Since the positive eigenvalues of $H$ are absent, it follows that the eigenvalues of $H$ are contained in a bounded set.  To prove Theorem \ref{th1.1} (b), it remains to show that  under the conditions (\ref{ass1})-(\ref{ass4}) with $ \rho_1'>2$ and $\rho_2>2$,  and that  zero is a regular point of $H$, the eigenvalues of $H$ can not accumulate to zero.
\\

 Let $\chi_1(x)^2 + \chi_2 (x)^2 =1$ be a partition of unity on $\bR^n$ with $\chi_1 \in C_0^\infty$  such that $\chi_1(x) = 1$ for $|x| \le R$ for some $R>0$. Let $-\Delta + \frac{q(\theta)}{r^2}$ denote its Friedrich's realization on the core $C_0^\infty(\bR^n \setminus\{0\}$. Set
\[
\tilde{R}_0 (z) = \chi_1 (- \Delta +1 -z)^{-1} \chi_1 + \chi_2 ( -\Delta + \frac{q(\theta)}{r^2} -z)^{-1} \chi_2,
\]
and
\[
K(z) = \tilde{R}_0 (z) (H-z)  - 1,
\]
for $z \in D'(0, \delta) = D(0, \delta) \setminus [0, \delta[$. One has
\begin{eqnarray*}
K(z) &= &    \chi_1 (- \Delta +1 -z)^{-1} ( [\Delta,\chi_1] +  (V-1) \chi_1)  \\
 &+ & \chi_2 ( -\Delta + \frac{q(\theta)}{r^2} -z)^{-1}  ([\Delta, \chi_2] + (V -  \frac{q(\theta)}{r^2}) \chi_2).
\end{eqnarray*}
Note that $(- \Delta +1 -z)^{-1} $ is holomorphic for $z$ near $0$ and the asymptotics of $( -\Delta + q(\theta)/r^2 -z)^{-1} $ can be computed by using Theorem \ref{thA.1} of Appendix A. In particular, (\ref{ass4}) ensures that the $\ln z$-term is absent.  One deduces that the limit \[
F_0 = \lim_{z \to 0, z\not\in \bR_+} \tilde{R}_0 (z)
 \]
 exists  in $\vL(-1,s; 1,-s)$ for any $s>1$.  Since $\rho_1'>2$ and $\rho_2>2$, $K(z)$ is a compact operator-valued function on $H^{1, -s}$, $1<s< \min\{ \rho'_1 /2, \rho_2/2\}$, holomorphic in $D'(0, \delta)$,  continuous up to the boundary  and
\[
 K(z) = K_0 + O(|z|^{\delta_0})
\]
in $\vL(1,-s; 1,-s)$ for some $\delta_0 >0$, where $K_0 =\lim_{z \to 0, z\not\in \bR_+} K(z)$ is a compact operator.  The assumption that zero is a regular point of $H$ implies that $-1$ is not an eigenvalue of $K_0$.  In fact, let $u \in H^{1,-s}$ for any $s>1$ such that $K_0 u =-u$. By the expression of $K_0$, one sees that $Hu = -HK_0 u  =0$ outside some compact. Therefore, 
 \[
H u =  \in  H^{-1, s} \mbox{ and }  F_0 H u=  (1+K_0) u =0.
 \]

Under the assumption (\ref{ass4}), using the decomposition of  $u$ in terms of the eigenfunctions of $-\Delta_{\bS^{n-1}} + q(\theta)$, one can show as in Theorem 3.1 of \cite{w2} that $u(x) = O(\w{ x}^{-\frac{n-2}{2} + \nu_0})$ for some $\nu_0>0$. Therefore, $u$ is in fact in  $ H^{1,-s'}$ for any $s'$ with  $ 1-\nu_0 < s'<1$.  On the other hand, one can show that the forms
\begin{eqnarray*}
q_1(v) &=& \lim_{  \lambda \to 0_-} \w{v, \chi_1 (- \Delta +1 -\lambda)^{-1} \chi_1  v} \quad \mbox{ and }  \\
q_2(v) &=& \lim_{  \lambda \to 0_-} \w{v,  \chi_2 ( -\Delta + \frac{q(\theta)}{r^2} -\lambda)^{-1} \chi_2  v}
\end{eqnarray*}
are positive  on $H^{-1,s}$ for any $s<1$ and $q_j(v)=0$ if and only if $\chi_j v=0$. Since $w = Hu$ verifies
\[
0 = \w{w, F_0 w} = q_1(w) + q_2(w)
\]
we deduce that $q_1(w)= q_2(w)=0$ and $\chi_1 w=0$ and $\chi_2 w=0$. This shows that $w=Hu=0$. The assumption that zero is a regular point of $H$ implies that $u=0$. This verifies that $-1$ is not an eigenvalue of $K_0$ on $H^{1,-s}$ for any $s>1$. In  particular, $( 1 + K_0)^{-1}$ existe in $\vL(1, -s; 1,-s)$. An argument of perturbation as used in the proof Theorem \ref{th1.1} (a) shows that
$(1 + K(z)$ is invertible for $z \in D'(0, \delta)$ with $\delta>0$ small enough. In particular, $  (1 + K(z))^{-1}$ has no poles there.  By Lemma \ref{lem2.1}, the eigenvalues of $H$ in $D'(0, \delta)$ as operator on  $L^2$ are the same as the poles of $(1+ K(z))^{-1}$  as operator on $H^{1,-s}$.  This shows that $H$ has no eigenvalues in $D'(0, \delta)$. \\

When $V$ is dissipative, Lemma \ref{lem2.2} shows that zero is a regular point of $H$ and the result of Theorem \ref{th1.1} (b) holds under the conditions (\ref{ass1})-(\ref{ass4}) with $\rho_1'>2$ and $\rho_2>2$. \ef
\bigskip

In the case where $-1$ is an eigenvalue of $K_0$,  one may try to use the Feshbach-Grushin formula
\begin{equation}
(1 + K(z) )^{-1} = E'(z) - (  1 - E'(z)K(z)) \Pi F(z)^{-1} \Pi ( 1 -  K(z) E'(z))
\end{equation}
where $F(z)$ is the Feshbach operator
\[
F(z) = \Pi(1+ K(z) - K(z) E'(z) K(z)) \Pi
\]
and to show that the zeros of the determinant of $ F(z)$ can not accumulate to zero. For real-valued potentials $V(x)$, this is indeed the case  because the inverse of $ F(z)$  can be explicitly calculated for $z$ near zero and $z \neq 0$ (\cite{jk,w2}). The same methods can be utilized when the imaginary part of the potential is sufficiently small. See Section 4.  In the general case,  new phenomena appear in threshold spectral analysis for complex-valued potentials, because although the kernels of $1+K_0$ and $H$ in $H^{1,-s}$ are the same,  the characteristic space of $K_0$ associated with the eigenvalue $-1$ is usually different from that of $H$ associated with the eigenvalue zero in $H^{1,-s}$.

\sect{Number of eigenvalues under perturbation}

 This Section is devoted to studying the number of eigenvalues under the perturbation by a small complex-valued potential. The main attention is payed to the perturbation of the zero eigenvalue and the zero resonance for dissipative operators. We first begin with the easy case where zero is a regular point of the selfadjoint operator $H_1$.
 \\

 Under  the conditions (\ref{ass1})-(\ref{ass4})  with $\rho_2>2$,
  $H_1= -\Delta_1 + V_1$ has only a finite  number of eigenvalues:
\begin{equation}
\nu_1 < \nu_2 <  \cdots < \nu_k  \le  0.
\end{equation}
The resolvent $R_1(z)$ of $H_1$ verifies for any $\delta >0$, $s>1/2$,
\begin{equation}
\| \w{x}^{-s} R_1(z)\w{x}^{-s}\| \le C_\delta |z|^{-\frac{1}{2}},
\end{equation}
for all $z \in \bC_-  $ with $|z| >\delta $ and $| z - \nu_j| > \delta$. If zero is a regular point of $H_1$, then  $\nu_k <0$ and one has for any $s>1$,
\begin{equation}
\| \w{x}^{-s} R_1(z)\w{x}^{-s}\| \le C,
\end{equation}
for all $z \in \bC_-  $ with $|z|  \le \delta $ for some $\delta >0$ small enough.
Here and in the following, $\| \cdot\|$ denotes the norm or the operator norm on $L^2$ with no possible confusion. Since $V_2$ is $-\Delta$-compact and is bounded by $C\w{x}^{-\rho_2}$ outside some compact and $\rho_2 >2$, we have
\begin{equation} \label{3.1}
\| |V_2|^{1/2} R_1(z)|V_2|^{1/2}\| \le C_\delta
\end{equation}
for all $z\in C_-$ with $| z-\nu_j| >\delta$, $j =1, \cdots, k$, if zero is a regular point of $H_1$. \\

\begin{prop}\label{prop4.1}
 Under the conditions (\ref{ass1})-(\ref{ass4}), let $H_1 = -\Delta + V_1$ and $H(\lambda) = H_1 - i \lambda V_2$, $\lambda \in \bR$, $V_2 \neq 0$.  Let  $N(\lambda)$ (resp., $N_1$) denotes the number of eigenvalues of $H(\lambda)$ (resp., $H_1$), counted according to their algebraic multiplicity.  Assume that zero is a regular point of $H_1$. Then there exists some $\lambda_0>0$ such that
\begin{equation}
N(\lambda) = N_1
\end{equation}
for $ |\lambda| < \lambda_0$.
\end{prop}
\pf
For each negative eigenvalue $\nu_j<0$ of $H_1$ with multiplicity $m_j$, a standard perturbation argument can be used to show that $\exists \delta_0>0$ such that $H(\lambda)$
has $m_j$ eigenvalues, counted according to their algebraic multiplicity,  in $D_-(\nu_j, \delta_0)$ if $\lambda >0$ is small enough.
Let $ \Omega = \bC_- \setminus (\cup_j D_-(\nu_j, \delta_0))$.  By (\ref{3.1}), there exists $\lambda_0>0$ is such that
\begin{equation} \label{3.2}
\lambda_0 \| |V_2|^{1/2} R_1(z)|V_2|^{1/2}\| <1,
\end{equation}
for all $z \in \Omega$.  Then, $H(\lambda)$ has no eigenvalues in $\Omega$ if $|\lambda| < \lambda_0$, because if $ u$ is an eigenfunction of $H(\lambda)$ associated with the eigenvalue $z_0  \in \Omega$,
\[
 H(\lambda) u = z_0 u,
\]
 then
$v = |V_2| ^{1/2}u \neq 0$.  $v$ is a non zero solution  of the equation
\[
v = -i \lambda |V_2|^{1/2} R_1(z_0) V_2^{1/2}v, \quad V_2^{1/2} = \mbox{ sign } V_2 | V_2|^{1/2}.
\]
This is impossible, because $\| \lambda |V_2|^{1/2} R_1(z_0) |V_2|^{1/2} \| <1$.
This proves that the total number of complex eigenvalues of $H(\lambda)$ is equal to the number of negative eigenvalues of $H_1$.
\ef

 Note that if $N_1=0$ and if zero is a regular point of $H_1$,  Kato's smooth perturbation argument can be directly used to show that $H(\lambda)$ is similar to $H_1$ for $\lambda $ small enough, which implies $N(\lambda)=0$. See \cite{ky,k}.\\

 Let us study now the perturbation of the threshold eigenvalue and resonance for dissipative Schr\"odinger operators. Assume that zero is an eigenvalue of $H_1 = -\Delta + V_1$ with multiplicity $k_0$ and a resonance of multiplicity $k$. $k_0$ or $k$ may eventually be  equal to $0$.   We want to show $H(\lambda) = H_1 - i \lambda V_2$ has $m=k_0 + k$ complex eigenvalues near $0$, counted according to their algebraic multiplicity.  The threshold eigenvalues and resonances are unstable under perturbation and may produce both eigenvalues or quantum resonances, if the potentials are dilation-analytic. Therefore,  it may be interesting to see why the zero eigenvalue and the zero resonance of $H_1$  will be turned into the eigenvalues of the non-selfadjoint Schr\"odinger operator $H$ as soon as a small dissipative part of potential appears. \\

\begin{theorem} \label{th4.1}  Assume (\ref{ass1})-(\ref{ass4})  with $\rho_1' >3$ and $\rho_2>3$ and (\ref{mu1}). Suppose that $V$ is dissipative: $V_2 \ge 0$ and $V_2 \neq 0$.

(a). If zero is an eigenvalue of multiplicity $m$, but not a resonance of $H_1$, then  there exists $\delta, \lambda_0>0$ such that for $0<\lambda <\lambda_0$, $H(\lambda)$ has exactly $m$ eigenvalues in $D_-(0, \delta)$.

(b). If zero is a resonance, but not an eigenvalue of $H_1$, suppose in addition  that the condition (\ref{ass5}) is satisfied. Then for $0<\lambda <\lambda_0$, $H(\lambda)$ has exactly one eigenvalue  in $D_-(0, \delta)$.
\end{theorem}
\pf
 Set
\begin{equation}
H_0 = -\Delta + \frac{q(\theta)}{r^2}, \quad  W_1(x) = V_1(x) -   \frac{q(\theta)}{r^2}
\end{equation}
so that $H_1 = H_0 + W_1$.  Let $R_j(z) = (H_j -z)^{-1}$, $j=0,1$.
 The low-energy asymptotic expansion of $R_0(z)$ can be explicitly calculated (\cite{w1}). One has
 \[
 R_0(z) = G_0 +  z_{\nu_1} G_1 + z G_2 +  O( |z|^{1+\epsilon})
 \]
 where  $G_j$ is continuous from $H^{-1, s}$ to $H^{1,-s}$ with $ s>j$. In particular, $G_1$ is a rank one operator given by
 \begin{equation}
 G_1 u = c_{0} \w{\eta_0, u}\eta_0, \quad u \in H^{-1,s}, s>2,
 \end{equation}
 with $\eta_0$ defined in Introduction and
 \begin{equation}
 c_0 = -\frac{e^{-i\pi \nu_1} \Gamma(1-\nu_1)}{\nu_1 2^{2\nu_1+1} \Gamma (1 + \nu_1)},  \mbox{ if } 0 < \nu_1<1; \quad c_0 = - \frac{1}{8}   \mbox{ if } \nu_1 =1.
 \end{equation}

 $z_{\nu_1}$ is defined by
 \[
 z_{\nu_1} = \left\{ \begin{array}{rl}
  e^{\nu_1 \ln z}, & \nu_1 \in ]0, 1[ \\
 z \ln z, & \nu_1 =1.
 \end{array}
 \right.
 \]
 with the branch of $\ln z$ chosen such that it is holomorphic on the slit complex plane $\bC\setminus \bR_+$  and
 $\lim_{\epsilon \to 0+} \ln (\lambda + i \epsilon) = \ln \lambda$ if $\lambda >0$. Although $W_1$ is not $-\Delta$-form compact, it is still $-\Delta$-form bounded for $n \ge 3$. By the comparaison with the decomposition
 \[
  H_1 =  H'_0 + W_1', \quad \mbox{ with } H_0' = \chi_1( -\Delta +1)\chi_1 + \chi_2 (-\Delta + \frac{q(\theta)}{r^2})\chi_2
 \]
where $\chi_1^2 + \chi_2^2 =1$ with $\chi_2(x)=0$ for $x$ near $0$, one can show  (\cite{j}) that
$1 + G_0W_1$ is a Fredholm operator and its kernel, $K$, in $H^{1,-s}$ coincides with that of $1+ G_0'W'_1$, where
\[
G_0' = \lim_{z\to 0, z\not\in \bR_+} (H_0'-z)^{-1}, \quad \mbox{ in } \vL(-1, s; 1,-s).
\]
 Using the positivity of $H_0$, one can show that the Hermitian form
 \[
 K\times K \ni (\varphi, \psi) \to \w{\varphi, -W_1 \psi} \in \bC
  \]
 is positive definite and  there exists a basis $\{\phi_1, \cdots, \phi_m\}$, $m=k_0 +k$,
 such that
 \[
 \w{\phi_i, - W_1 \phi_j } =\delta_{ij}.
 \]
The assumption (\ref{mu1}) implies that the multiplicity of the zero resonance is at most one. If zero is a resonance, we may assume that $\phi_1$ is a resonant state and the others are eigenfunctions of $H_1$. Let $ Q : H^{1, -s} \to H^{1,-s}$, $Qf =\sum_{j=1}^m \w{-W_1\phi_j, f} \phi_j$. $Q$ is a projection from $H^{1, -s}$ onto $K$. Set $Q' =1-Q$. One can show that range of $Q'$ is closed and is equal to the range of $1 + G_0 W_1$ in $H^{1,-s}$. Then the Fredholm theory  shows that $(Q'(1 +  G_0 W_1)Q')^{-1}Q'$ exists and is continuous on $H^{1, -s}$. See \cite{j,w2}. \\

 Let $R(z, \lambda) = (H(\lambda) -z)^{-1}$, $z \not \in  \sigma (H(\lambda)$.  One has the resolvent equations
\[
R(z, \lambda ) = (1 -i \lambda R_1(z) V_2)^{-1}R_1(z) = (1 + R_0(z)(W_1- i \lambda V_2))^{-1}R_0(z).
\]
 By Lemma \ref{lem2.1}, the eigenvalues of $H(\lambda)$ in $D_-(0, \delta)$ are the same as the set of $z$ such that $Hu = zu$ has a solution $u$ in $H^{1, -s}$ for some $s>0$.   Since $R_0(z)$ is holomorphic in $\bC_-$, a point $z_0 \in \bC_-$ is an eigenvalue of $H(\lambda)$ if and only if it is a pole of
\[
z \to  (1 + R_0(z) (W_1 - i \lambda V_2))^{-1}
\]
and their multiplicity is the same. Let
\[
W(z, \lambda) = 1 + R_0(z) (W_1 - i \lambda V_2).
\]
For $\lambda >0$  small enough, a perturbation argument shows that
\[
E'(z,  \lambda) := (Q'W(z, \lambda)Q')^{-1} Q'
 \]
 exists on $H^{1, -s}$ for $1<s< 1 +\epsilon$  and  $z\in D_-(0,, \delta)$,
 $\| E'(z, \lambda)\|_{\vL(1,-s; 1,-s)}$ is uniformly bounded for $z \in D_-(0, \delta)$ and $\lambda$.
This allows us to show that the Grushin problem
\[
\vW(z, \lambda)=\left( \begin{array}{cc}
                  W(z, \lambda) & T\\[.1in]
                   S & 0
                \end{array}
                \right) \;  :\; H^{1,-s}\times \bC^m \to H^{1,-s}\times \bC^m,
\]
where $s>{1}$,   $T$ and $S$  are given by
\begin{eqnarray*}
T   c &=& \sum_{j=1}^m c_j\phi_j, \quad  c = (c_1,
                  \cdots, c_m)\in \bC^m, \\
Sf &= &(\w{-W_1\phi_1, f},  \cdots, \w{-W_1\phi_m, f}) \in \bC^m,\quad
f\in H^{1,-s}.
\end{eqnarray*}
The inverse of $\vW(z, \lambda)$ is given by
\[
\vW(z, \lambda)^{-1}=\left(
                           \begin{array}{cc}
              E(z) & E_+(z)
                            \\[.1in]
             E_-(z)  & E_{+-}(z)
                           \end{array}
                          \right),
\]
where
\begin{eqnarray}
E(z) & = & E'(z), \\
 E_+(z) &=& T - E'(z) W(z)T, \\
E_-(z )  &=&  S -SW(z) E'(z),\\
  E_{+-}(z)& =& - SW(z)T + S W(z) E'(z) W(z)T.
\end{eqnarray}
Here to simplify typesetting, the indication of dependance on $\lambda$ is omitted.
It follows that the inverse of $W(z)$ is given by
\begin{equation} \label{l}
W(z)^{-1} = E(z) - E_+(z) E_{+-}(z)^{-1} E_-(z) \mbox{ on } H^{1, -s}.
\end{equation}
Since $E(z),  E_\pm(z)$ and $E_{+ -}(z)$ are holomorphic and uniformly bounded  as operators on $H^{1,-s}$ for $z \in D_-(0, \delta)$ and $|\lambda| <\delta$,   $z_0$ is a pole of  $W(z)^{-1}$ in $D_-(0, \delta)$  if and only if
\[
F(z_0, \lambda) := \det E_{+ -}(z_0) =0
\]
and their multiplicities are the same (\cite{hs}). As seen before,  by Lemma \ref{lem2.1}, the poles of $ z  \to W(z)^{-1}$ in $\bC_-$ as operator on $H^{1,-s}$ are the eigenvalues of $H$.
\\

We are now led to prove that for $0 <\lambda \le \lambda_0$ small enough,  $F(z, \lambda)$ has $m$ zeros in $D_-(0, \delta)$.
Since $\phi_j \in L^{2,-s}$ for any $s>1$, under the condition $\rho_1' >3$ and $\rho_2>3$ and that zero is a resonance or an eigenvalue of $H_1$,  one can calculate  the asymptotics of the matrix
\begin{equation}
E_{+-}(z)  = \left(\w{-W_1\phi_j,  (W(z) -  W(z) E'(z) W(z))\phi_k}\right)_{1\le j, k \le m},
\end{equation}
up to an error $O(|z|^{1+\epsilon})$. This computation is known in the case $\lambda =0$ (see Proposition 4.4, \cite{w2}). In the case $\lambda \neq 0$, the calculation is similar. We give only the result and omit the details.
\begin{equation}
(E_{+-}(z))_{j,k} = -i \lambda v_{jk} + z_{\nu_1} a_{jk} + z (b_{jk} + r_{jk}) + O(|z|^{1+\epsilon})
\end{equation}
where
\begin{eqnarray}
v_{jk} & = & \w{\phi_j, V_2 \phi_k} \\
a_{jk} & = & -c_0 (|c_1|^2 \delta_{1j}\delta_{1k}  -i \lambda \overline{c}_1 c_k' \delta_{1k} )\quad  \mbox{ with } \\
c_1 & = & \w{ W_1 \eta_0, \phi_1}, \quad c_j' = \w{V_2 \eta_0, \phi_j},  \\
b_{jk} &=&  \w{-W_1 \phi_j,  G_2 W_1 \phi_k} + i \lambda \w{W_1\phi_j, G_2 V_2 \phi_k} \\
r_{jk} & =& \w{-W_1\phi_j, G_1 (W_1 - i \lambda V_2) E'(0) G_1 (W_1 - i \lambda V_2) \phi_k}.
\end{eqnarray}
Here $\delta_{1j} =1$ or  $0 $  according to $ j=1$  or  not.
 In the case $0$ is not a resonance, $c_1  =0$.  Remark also that if $\phi_j$ and $\phi_k$ are both eigenfunctions of $H_1$,
\begin{equation}
\w{W_1\phi_j, G_2 W_1 \phi_k} = \w{\phi_j, \phi_k}
\end{equation}
and if $\phi_j$ is an eigenfunction of $H_1$,
\begin{equation}
r_{jk} = 0, \quad k=1, \cdots, m.
\end{equation}
The matrice $(v_{jk})$ is positive definite.  In fact, it is clearly positive since $V_2 \ge 0$. If $0$ is an eigenvalue of this matrix, we can take an associated   eigenvector $d =(d_1, \cdots, d_m) \in \bC^m \setminus\{0\}$. Let $\psi = \sum_{j} d_j \phi_j$. One has $V_2\psi =0$, which implies
$\psi =0$ on a non trivial open set. Since $H_1\psi =0$, the unique continuation theorem implies $\psi =0$ on $\bR^3$. This leads to a contradiction with the fact that $\phi_1, \cdots, \phi_m$  are linearly independent. \\

Consider first the case zero is an eigenvalue but not a resonance of $H_1$.   One has $a_{jk} =0$ and
\begin{equation}
(E_{+-}(z))_{j,k} = -i \lambda v_{jk}   + z ( - \w{\phi_j, \phi_k} +  i \lambda \w{W_1 \phi_j, G_2 V_2 \phi_k} )  + O(|z|^{1+\epsilon})
\end{equation}
Since the matrices $(\lambda v_{jk}) $ and $(\w{\phi_j, \phi_k} )$ are positive definite,
$F_0(z, \lambda) = \det ( i (\lambda v_{jk})   + z  (\w{\phi_j, \phi_k} )) $ has m zeros of the form $ z = - i \lambda \sigma_j$, $\sigma_j>0$. Let $ - i \lambda \sigma$ be one of the zeros of $F_0(z, \lambda)$ with order $k$. For some appropriate $c>0$ such that the distance from zeros of $z \to F_0(z, \lambda)$ to the circle $\partial D( -i \lambda \sigma, c \lambda)$ is at least $c' \lambda$ for some $c'>0$, one has
\[
|F_0(z, \lambda) |\ge C_1 \lambda^m,  \quad |F(z, \lambda) - F_0(z, \lambda)| \le C_2 \lambda^{m+\epsilon}
\]
for $|z + i \lambda \sigma | = c \lambda$.
For $\lambda >0$ small, we can apply the Rouché's theorem to $F(z, \lambda) $ to conclude that $F(z, \lambda)$ has also $k$ zeros in the disk $D(- i \lambda \sigma, c\lambda)$.  Therefore, the total number of zeros of $F(z, \lambda)$ in $D_-(0, \delta)$ are at least $m$.
If $z_0$ is a zero of $F(z, \lambda)$, the asymptotic expansion of $F(z, \lambda)$ in $z$ shows that $ \varsigma = \lim_{\lambda \to 0} z_0/\lambda $ exists and $\lambda \varsigma$ is a zero of $F_0(z, \lambda)$.  This allows to conclude that $F(z, \lambda)$ has in all $m$ zeros in $D_-(0, \delta)$ for $\delta $ sufficiently small and $0 <\lambda <\lambda_0$ with $\lambda_0>0$ small enough.
\\

 Assume now that zero is a resonance, but not an eigenvalue of $H_1$. $E_{+-}(z)$ is a scalar function, holomorphic in $z \in D_-(0, \delta)$ and
\begin{equation}
E_{+-}(z) = - i \lambda v_{11} - c_0  z_{\nu_1} (|c_1|^2   -i \lambda \overline{c}_1 c_1' ) + O(|z|).
\end{equation}
Assume (\ref{ass5}). The root of the equation
\begin{equation} \label{p}
i \lambda v_{11} +  c_0  z_{\nu_1} (|c_1|^2   -i \lambda \overline{c}_1 c_1' ) =0.
\end{equation}
can be explicitly calculated. Set
\[
r = \frac{ \lambda v_{11}}{ | c_0 | \; |( |c_1|^2   -i \lambda \overline{c}_1 c_1' )| }.
\]
Since  $\overline{c}_1 c_1'  <0$ by the assumption (\ref{ass5}),  the argument, $\varphi$, of
$-\lambda  \overline{c}_1 c_1'  + i |c_1|^2$  verifies
\begin{equation}
0 < \varphi < \frac{\pi}{2}, \quad \lim_{\lambda \to 0 _+} \varphi = \frac{\pi}{2}.
\end{equation}
For $\nu_1 \in ]0, 1[$, $z_{\nu_1} = z^{\nu_1}$ and  $c_0= -\frac{  e^{-i \pi \nu} \Gamma(1-\nu)}{ \nu 2^{2\nu+1}  \Gamma ( 1 + \nu)}$. With the above definitions of $r$ and $\varphi$, (\ref{p}) can be rewritten as
\begin{equation} \label{p1}
  z_{\nu_1} = - \frac{  i \lambda v_{11}}{ c_0  (|c_1|^2   -i \lambda \overline{c}_1 c_1' ) } = r e^{i (\pi  \nu_1 + \varphi)}.
\end{equation}
When $   \nu_1 \in[1/2, 1[$,  (\ref{p1}) has a unique solution $z_0 = \rho e^{i\theta}$ in $D_-(0, \delta)$ given by
\begin{equation}
\rho = r^{1/\nu_1} \quad   \mbox{ and }  \quad  \theta =  \pi  +  \nu_1^{-1}  \varphi.
\end{equation}
 When $\nu_1=1$, $z_{\nu_1} = z \ln z$ and $c_0 = -\frac 1 8$. $ z = \rho e^{i \theta}$ is a solution of (\ref{p}) if
 \begin{equation} \label{p2}
 \rho e^{i \theta}  (\ln \rho + i\theta)  = r e^{i \varphi}
 \end{equation}
 Remark that $r \to 0$ as $\lambda \to 0$ and $\arg  (\ln \rho + i\theta) \to \pi_-$ as $\rho \to 0_+$.  Set  $\ln z = \tau e^{i (-\pi -\sigma)}$ with $\tau = |\ln z |$. For $\lambda >0$ small enough, one can check that the system
 \[
 \left\{ \begin{array}{ccl}
  \tau e^{-\tau \cos \sigma} & = & r, \\
  -\sigma + \tau \sin \sigma & = & \varphi + \pi
 \end{array}
 \right.
 \]
 has a unique solution $(\tau, \sigma)$ such that $\tau \to \infty$ and  $\sigma \to 0_+$ as $\lambda \to 0$. This shows  that (\ref{p2}) has a unique solution $z_0= \rho e^{i\theta}$
 in $D_-(0, \delta)$ given by
 \[
 \rho = | e^{\tau e^{i (-\pi -\sigma)}}| =  e^{-\tau \cos \sigma} \quad \mbox{  and } \quad \theta = \pi + \varphi-\sigma.
 \]

  Using the Rouché's theorem, one can show as before that $E_{+-}(z, \lambda)$ has just one zero in $D_-(0, \delta)$ which is located  inside the small disk $D(z_0, g(\lambda))$, where $g(\lambda) = a \lambda^{1/\nu_1}$ if $\nu_1 \in ]0, 1[$ and $g(\lambda) = a \lambda \w{\ln \lambda}^{-1}$ if $\nu_1=1$ and $a>0$ is an appropriate constant. This proves that $H(\lambda)$ has exactly one eigenvalue in $D_-(0, \delta)$ and it is simple.
 \ef

Theorem \ref{th1.2} follows  from Theorem \ref{th4.1} and the argument used Proposition \ref{prop4.1} outside a small neighborhood of zero. \\

 The case that zero is both an eigenvalue and a resonance can in principle be analyzed in a similar way, using the asymptotic expansion of $E_{+-}(z, \lambda)$ given above under stronger decay assumption $\rho_1'>4$ and $\rho_2>4$. But the factorization and the evaluation of  zeros of the determinant raise some difficulties when $m$ is arbitrary. We do not go further here.  Without the assumption (\ref{mu1}),   the zero resonance of $H_1$ may appear of arbitrary multiplicity. The matrix $E_{+-}(z, \lambda)$ can be still calculated, but its analysis is more complicated. \\

 \begin{rmk} An example for which  zero is not an eigenvalue, but a resonance of $H_1 = -\Delta + V_1$ can be constructed as follows.   Let $n=3$ or $4$ and let $v_1 $ be a real-valued function satisfying (\ref{ass1})  with $\rho_1 > 3$ and $v_1 \le 0$, $v_1 \neq 0$. Let $H_1(\beta) = - \Delta + \beta v_1$. Then one can show that  there exists a critical constant $\beta_0>0$ such that $ H_1(\beta_0)  \ge 0$ and $ H_1(\beta)  $ has at least one negative eigenvalue for any $\beta >\beta_0$. Then zero is a resonance but not an eigenvalue of $ H_1(\beta_0) $. Take $V_1 =\beta_0 v_1$.  A resonant state of $H_1 = -\Delta + \beta_0 v_1$ can be constructed as some weak limit of the fundamental state of $ H_1(\beta)$ as  $\beta \to \beta_{0+}$ and it does not change sign. The condition (\ref{ass5}) is then satisfied  for any $V_2 \in C_0^\infty(\bR^n)$ with $V_2 \ge 0$ and $V_2 \neq 0$. In this example, the number of eigenvalues of $H_1$ is zero, while that of $H_1 - i \lambda V_2$ is one  for any $\lambda>0$ small enough, by Theorem \ref{th1.2}. \\
\end{rmk}

\appendix{ \sect{Low-energy resolvent expansion on conical manifolds}  \label{app} }

In this appendix, we recall in a concise way  the result of \cite{w1} used in Section 4 on
low-energy  resolvent expansion of  the model operator and at the same time correct some sign error. Consider the operator
\begin{equation}
P_0 = - \Delta_{g} + \frac{q(\theta)}{r^2}
\end{equation}
on a conical manifold  $M =\bR_+ \times \Bsigma$ equipped with a Riemannian metric $g$, where $\Bsigma$ is an
$(n-1)$-dimensional compact manifold, $n\ge 2$. Here $(r, \theta)\in \bR_+ \times\Bsigma$, $q(\theta)$ is a real
continuous function and the metric  $g$ is of the form
\[
g = dr^2 + r^2 h
\]
with $h$ a Riemannian metric on $\Bsigma$ independent of $r$. If
$\Bsigma$ is of boundary, the Dirichlet condition is used for
$P_0$.  We still denote by  $P_0$  its Friedrich's realization  with the core $C_0^\infty (\bR^n \setminus \{0\})$. Let  $\Delta_h$ denote Laplace-Beltrami operator on
$\Bsigma $.  Assume
\begin{equation}
-\Delta_h + q(\theta) \ge -\frac{(n-2)^2}{4},  \quad  \mbox{ on } L^2(\Bsigma).
\end{equation}
Put
\begin{equation}\label{nu}
\sigma_{\infty} =\left\{ \nu; \; \nu =\sqrt{ \lambda +
\frac{(n-2)^2}{4}}, \lambda \in \sigma (-\Delta_h + q) \right\}.
\end{equation}
Denote
\[
\sigma_k  = \sigma_{\infty} \cap [0,k ], \quad \quad k \in \bN.
\]
Let $\pi_\nu$ denote the orthogonal projection in $L^2(M)$ onto
the subspace spanned by the eigenfunction  of $-\Delta_h + q$
associated with the eigenvalue $\lambda_\nu = \nu^2 - \frac{(n-2)^2}{4}$.
Define for $\nu \in \sigma_{\infty}$
\[
z_\nu =\left\{ \begin{array}{ccl}
z^{\nu'}, &  \quad & \mbox{ if } \nu\not \in\bN \\
  z \ln z, & & \mbox{ if } \nu \in\bN^*,
 \end{array}
 \right.
\]
where $\nu' =\nu -[\nu]$.
Let $\sigma_N =\sigma_{\infty}\cap [ 0, N]$. For $\nu>0$,  let
$[\nu]_-$ be the largest integer strictly less than $\nu$. When
$\nu=0$, set $[\nu]_-=0$. Define $\delta_\nu$ by $\delta_\nu = 1$,
if $\nu \in \bN\cap \sigma_\infty$; $0$, otherwise.  \\

\begin{theorem} \label{thA.1}
 Let $R_0(z) = (P_0-z)^{-1}$ for $z\not \in \bR_+$. The following asymptotic
expansion holds for $z$ near $0$ with $ z
\not\in \bR_+$.
\begin{equation}\label{libre}
R_0(z)  =   \delta_0 \ln z \; G_{0,0} +  \sum_{j=0}^{N} z^j F_j
+ \sum_{\nu \in \sigma_N} z_{\nu}\sum_{j=[\nu]_-}^{N-1} z^j
G_{\nu,j+ \delta_\nu }\pi_\nu  + O( |z|^{N +\epsilon}),
\end{equation}
 in $ \vL(-1,s; 1, -s), \;  s > 2N + 1$.  Each term in the above expansion  can be explicitly calculated. In particular, $F_j \in \vL(-1,s; 1, -s), \;  s > 2j + 1$ and $G_{\nu, j}$, $j \ge [\nu]_-$,  is of finite rank with its  Schwartz kernel on $L^2(\bR_+; r^{n-1} dr)$ given by
\begin{equation} \label{a.5}
G_{\nu, j}(r, \tau)  =   \frac{  (-1)^{j+1 -[\nu]}  e^{-i \nu' \pi}\Gamma (1-\nu') \; \; (r\tau)^{-\frac{ n-2}{2} + \nu' +j} }{  2^{2 \nu +1} \pi^{\frac 1 2}  (j-[\nu])! \; \Gamma (\frac  1 2 +\nu) \nu' (\nu'+1)\cdots (\nu'+j) }  P_{\nu,j-[\nu]}(\rho)
\end{equation}
for $\nu \not\in \bN$, $ j \ge [\nu]$  and $\nu'  = \nu -[\nu]$; and
\begin{equation} \label{a.6}
G_{\nu, j}(r, \tau)  =   \frac{  (-1)^{j+l +1}  }{ \pi^{\frac 1 2} 2^{2 l +1} j! \; (j-l)! \;  \Gamma (l + \frac 1 2)} (r\tau)^{-\frac{ n-2}{2}  +j}  P_{l,j- l}(\rho)
\end{equation}
for $\nu = l \in \bN$ and $j \ge l$. $P_{\nu, k}(\rho)$ is a polynomial of degree $k$ in $\rho$:
\begin{equation}
P_{\nu, k}(\rho) = \int_{-1}^1 ( \rho + \frac \theta 2)^k ( 1 -\theta^2)^{\nu-1/2} d\theta, \quad \rho = \frac{r^2 + \tau^2}{4r \tau}.
\end{equation}
\end{theorem}

 Note that
\begin{equation}
P_{\nu,0}(\rho) = \frac{ \Gamma(1/2) \Gamma ( 1/2 + \nu)}{ \Gamma (1 + \nu)}, \quad  P_{ \nu,1}(\rho) =
\frac{ \Gamma(1/2) \Gamma ( 1/2 + \nu)}{ \Gamma (1 + \nu)} \rho.
\end{equation}
One has for $\nu \in  [0, 1]$
\begin{equation}
G_{\nu, [\nu]}(r, \tau) = \gamma_\nu (r\tau)^{-\frac{n-2}{2} + \nu},
\end{equation}
with
\begin{equation}
\gamma_\nu  =  -\frac{  e^{-i \pi \nu} \Gamma(1-\nu)}{ \nu 2^{2\nu+1}  \Gamma ( 1 + \nu)} \; \mbox{ for }  \;  \nu \in ]0, 1[, \quad  \gamma_0 = - \frac{1}{2}  \quad \mbox{ and } \quad
\gamma_1   =  -\frac{1}{8}.
\end{equation}

The expansion of $R_0(z)$ is obtained by decomposing $R_0(z)$ into
\[
 R_0(z) =  \sum_{\nu \in \sigma_\infty} (Q_\nu -z)^{-1} \pi_\nu
\]

with
\[
 Q_{\nu} = - \frac{d^2}{dr^2} -  \frac{n-1}{r}  \frac{d}{dr} +  \frac{\nu^2-
\frac{(n-2)^2}{4}}{r^2}, \quad \mbox{ on }
L^2(\bR_+;  r^{n-1} dr).
\]
Since the Schwartz kernel of $e^{-it Q_\nu}$ is given by
\begin{equation}
  \frac{1}{2it} (r\tau )^{ -\frac{n-2}{2}} e^{ -
\frac{r^2 + \tau^2}{4it}  -i \frac{\pi \nu}{2}} J_\nu(
\frac{r\tau }{2t}), \quad t\in \bR,
\end{equation}
where $J_\nu(\cdot)$ is the Bessel function of the first
kind of order $\nu$ and
\[
(Q_\nu -z)^{-1} = i \int_0^\infty e^{-it (Q_\nu-z)} \; dt
\]
for $\Im z>0$,  the Schwartz kernel of $(Q_\nu -z)^{-1}$ is
\begin{eqnarray}\label{Knu}
K_\nu(r, \tau  ; z) & = &  (r\tau )^{ -\frac{n-2}{2}}\int_0^\infty e^{ -
\frac{r^2 + \tau^2}{4it} + izt -i \frac{\pi \nu}{2}} J_\nu(
\frac{r\tau }{2t}) \;  \frac{dt}{2t}.
\end{eqnarray}
for $\Im z > 0$. The formula (2.6) in \cite{w1} for this kernel contains a wrong sign.  The coefficients in (\ref{a.5}) and (\ref{a.6})  are obtained from the constants given in Section 2 and Appendix A of \cite{w1}, in taking into account this sign correction.   Note that under the assumption (\ref{ass5}),  the continuity of $F_j$ and the remainder estimate can be improved. See Remark 2.4 in \cite{w1}.

\end{document}